\documentclass{amsart}

\usepackage{amsfonts}
\usepackage{amsmath}
\usepackage{amsthm}
\usepackage{amssymb}
\usepackage[english]{babel}
\usepackage{graphics}
\usepackage{latexsym}
\usepackage{longtable} 
\usepackage{mathrsfs} 
\usepackage{graphicx}
\usepackage{wrapfig} 
\usepackage[pdftex,colorlinks=true]{hyperref}

\usepackage{caption}
\usepackage{subcaption}

\newtheorem{theorem}{Theorem}
\newtheorem{corollary}[theorem]{Corollary}
\newtheorem{lemma}[theorem]{Lemma}

\newtheorem{claim}[theorem]{Claim}
\newtheorem{example}[theorem]{Example}

\theoremstyle{definition}
\newtheorem{definition}[theorem]{Definition}
\newtheorem{remark}[theorem]{Remark}

\newcommand{\mO}{\mathcal{O}}

\newcommand{\mD}{\mathcal{D}}

\newcommand{\E}{\mathrm{E}}
\newcommand{\D}{\mathrm{D}}

\newcommand{\R}{\mathbb{R}}

\newcommand{\noi}{\noindent}
\newcommand{\ms}{\medskip}

\newcommand{\De}{\Delta}

\newcommand{\la}{\lambda}

\newcommand{\Om}{\Omega}

\newcommand{\larrow}{\longrightarrow}

\newcommand{\sub}{\subseteq}

\newcommand{\bt}{\begin{theorem}}\newcommand{\et}{\end{theorem}}
\newcommand{\bd}{\begin{definition}}\newcommand{\ed}{\end{definition}}
\newcommand{\bl}{\begin{lemma}}\newcommand{\el}{\end{lemma}}
\newcommand{\beq}{\begin{equation}}\newcommand{\eeq}{\end{equation}}
\newcommand{\bc}{\begin{claim}}\newcommand{\ec}{\end{claim}}
\newcommand{\bex}{\begin{example}}\newcommand{\eex}{\end{example}}
\newcommand{\bcor}{\begin{corollary}}\newcommand{\ecor}{\end{corollary}}
\newcommand{\bp}{\begin{proof}}\newcommand{\ep}{\end{proof}}

\newcommand{\BPT}{\medskip \noindent \textbf{Proof of Theorem} }

\numberwithin{equation}{section}

\begin{document}


\title[Explicit $\infty$--Harmonic functions in high dimensions]{Explicit $\infty$-Harmonic functions in high dimensions}

\author{Birzhan Ayanbayev}
\address{Department of Mathematics and Statistics, University of Reading, Whiteknights, PO Box 220, Reading RG6 6AX, Berkshire, UK}



\email{b.ayanbayev@pgr.reading.ac.uk}



\date{}



\begin{abstract} 
The aim of this work is to derive new explicit solutions to the $\infty$-Laplace equation, the fundamental PDE arising in Calculus of Variations in the space $L^\infty$. These solutions obey certain symmetry conditions and are derived in arbitrary dimensions, containing as particular sub-cases the already known classes two-dimensional infinity-harmonic functions.
\end{abstract}

\maketitle


\section{Introduction} \label{section1}
Let $\Om \sub \R^n$ be an open set and $u\in C^2(\Om)$ a continuous twice differential function. In this paper we study the existence of solutions to the PDE
\beq\label{1}
\De_\infty u\, :=  \sum^n_{i,j=1}\D_i u \, \D_j u \, \D^2_{ij}u \, =\, 0
\eeq
of the form
\[
u(x) = \displaystyle\prod^n_{i=1}f_i(x_i),
\]
where $f_i$ are possibly non-linear for $1\leq i \leq n$, and $x = (x_1,...,x_n)^\top$, $x\in\Om$. Solutions of this form are called $separated$  $\infty$-$harmonic functions$. In the above $\D_i \smash \equiv \frac{\partial}{\partial x_i}$ and $\D^2_{ij} \equiv \smash{\frac{\partial^2}{\partial x_i\partial x_j}}$. 
The equation \ref{1} is called $\infty$-Laplacian (being a special case of the so-called more general the Aronsson equation) and it arises in Calculus of Variations in $L^\infty$ as the analogue of $Euler$-$Lagrange$ equation of the functional
\[
\E_\infty(u,\mO):=||\D u||_{L^\infty(\mO)}, \ \ \ \ \mO \Subset \Om,\ \ \ \ u\in W^{1,\infty}_\text{loc}(\Om,\R).
\]
These objects first arose in the work of G. Aronsson in the 1960s (see \cite{A1},\cite{A2}) and nowadays this is an active field of research for vectorial case $N\ge 2$ for $u\in W^{1,\infty}_\text{loc}(\Om,\R^N)$ which has begun much more recently in 2010s (see e.g. \cite{K1}). Since then, the field is developed enormously by N. Katzourakis in the series of papers (\cite{K3}-\cite{K11}) and also in collaboration with the author, Abugirda, Croce, Manfredi, Moser, Parini, Pisante and Pryer (\cite{AyK}, \cite{AK}, \cite{CKP}, \cite{KM}, \cite{KMo}, \cite{KPa}, \cite{KP1} - \cite{KP3}). A standard difficulty of \ref{1} is that it is nondivergence and since in general smooth solutions do not exist, the definition of generalised solutions is an issue. To this end, the theory of viscosity solutions of Crandall-Ishii-Lions is utilised (see e.g.\ \cite{K2}). 

In this paper all the $separated$ $\infty$-$harmonic functions$ are found for $n=2$ in polar coordinates, for $n=3$ in spherical coordinates and for all $n\geq 2$ in cartesian coordinates. Some of these new solutions derived herein coincide with previously known classes of solutions. For instance, the well-known G. Aronsson's solution $u(x,y) = |x|^\frac{4}{3} - |y|^\frac{4}{3}$ which has a $C^{1,1/3}$ regularity, described in Remark \ref{Remark 1}. Also M.-F. Bidaut-Veron, M. Garcia Huidobro and L. Veron have found a solution (\cite{VHV}) which is coincides with first two solutions of the theorem \ref{theorem1} and I.L. Freire, A. C. Faleiros have found a solutions of \ref{1} in \cite{FF}, but only one non-trivial of their solutions coincides with a particular case of Theorem \ref{theorem2} when $A = 1$. There may exist other additional solutions but this topic is not discussed herein.

The main results of this paper are contained in the following theorems.

\begin{theorem}[Separated two-dimensional $\infty$-Harmonic functions in polar coordinates] \label{theorem1} Let $u :\Om \sub\R^2 \larrow \R$ be a $separated$  $\infty$-$harmonic function$ of the $\infty$-Laplace equation in polar coordinates
\beq \label{1.1}
u^2_r u_{rr} + \frac{2}{r^2}u_r u_\theta u_{r\theta} + \frac{1}{r^4} u^2_\theta u_{\theta\theta} - \frac{1}{r^3}u_r u^2_\theta = 0
\eeq
of the form $u(r,\theta) = f(r)g(\theta)$, where $f, g$ are non-linear. Then, one of the following holds: either

\ \ \ $\mathrm{(i)}$  $|f(r)|=r^A$ and $|g(\theta)| = e^{B\theta}$, where $A$ and $B$ any constants, such that 
\[
A^2-A+B^2=0
\]

\noi or 

\ \ \ $\mathrm{(ii)}$  $|f(r)|= r^A$ and $|g(\theta)| = |g(\theta_0)| e^{\int^{\theta}_{\theta_0} G(t)d t}$, where $G(t)$ satisfies the following
\[
    t + c = 
\begin{cases}
    - \arctan\frac{G(t)}{A} +  \frac{A-1} {\sqrt{A^2-A}}\arctan\frac{G(t)}{\sqrt{A^2-A}}, & \text{if } A^2-A > 0\\
    \frac{1}{G(t)}, & \text{if } A = 0\\
    -\arctan{G(t)}, & \text{if } A = 1\\
    - \arctan\frac{G(t)}{A} +  \frac{A-1}{2 \sqrt{A-A^2}}\ln{\Big|\frac{G(t)-\sqrt{A-A^2}}{G(t)+\sqrt{A-A^2}}\Big|},              & \text{if } A^2-A < 0\\
\end{cases}
\]
\noi or

\ \ \  $\mathrm{(iii)}$  $|g(\theta)| = e^{B\theta}$ and $|f(r)| =  |f(r_0)| e^{\int^r_{r_0} \frac{\Phi(t)}{t} d t}$, where $\Phi(t)$ satisfies the following
\[
    \ln |t| + c = 
\begin{cases}
   \frac{1}{2}\ln\Big|\frac{\Phi^2(t)+B^2}{\Phi^2(t) - \Phi(t) + B^2}\Big| - \frac{1}{2}\frac{1}{\sqrt{B^2-\frac{1}{4}}}\arctan\frac{\Phi(t)-\frac{1}{2}}{\sqrt{B^2-\frac{1}{4}}}, & \text{if } B^2-\frac{1}{4} > 0\\
   \frac{1}{2}\ln\Big|\frac{\Phi^2(t)+B^2}{\Phi^2(t) - \Phi(t) + B^2}\Big| + \frac{1}{2} \frac{1}{\Phi(t)-\frac{1}{2}}, & \text{if } B^2-\frac{1}{4} = 0\\
   \frac{1}{2}\ln\Big|\frac{\Phi^2(t)+B^2}{\Phi^2(t) - \Phi(t) + B^2}\Big| - \frac{1}{4\sqrt{\frac{1}{4}-B^2}}\ln\Big|\frac{\Phi(t)-\frac{1}{2}-\sqrt{\frac{1}{4} - B^2}}{\Phi(t)-\frac{1}{2}+\sqrt{\frac{1}{4} - B^2}}\Big|, & \text{if } B^2-\frac{1}{4} < 0,\\
\end{cases}
\]
where c is any constant.
\end{theorem}

\begin{theorem}[Separated three-dimensional $\infty$-Harmonic functions in spherical coordinates] \label{theorem4}
Let $u:\Om\sub \R^3\larrow\R$ be a $separated$  $\infty$-$harmonic function$ of the $\infty$-Laplace equation in spherical coordinates
\beq\label{1.16}
\begin{split}
& u^2_r u_{rr} + \frac{2}{r^2 \sin^2\alpha}u_r u_\theta u_{r\theta} + \frac{1}{r^4\sin^4\alpha} u^2_\theta u_{\theta\theta} - \frac{1}{r^3\sin^2\alpha}u_r u^2_\theta + \\
& \frac{2}{r^2}u_r u_{\alpha} u_{r \alpha} + \frac{1}{r^4}u^2_{\alpha} u_{\alpha\alpha} - \frac{1}{r^3}u_r u^2_{\alpha}+\frac{2}{r^4 \sin^2 \alpha}u_{\theta}u_{\alpha}u_{\theta\alpha}-\frac{\cos\alpha}{r^4\sin^3\alpha}u_{\alpha}u^2_\theta= 0\\
\end{split}
\eeq
of the form $u(r,\theta,\alpha) = f(r) g(\theta) h(\alpha)$, where $f, g$ and $h$ are non-linear. Then, one of the following holds: either

\ \ $\mathrm{(i)}$  $|f(r)|=r^A$ and $|g(\theta)| = e^{B\theta}$, where $A$ and $B$ any constants, $h(\alpha) = h(\alpha_0)e^{\int^{\alpha}_{\alpha_0} H(t) d t}$, where $H(t)$ satisfies the following
\[
\Big|A^2-A+\frac{B^2}{\sin^2 t}+H^2(t)\Big| = c(A,B,t_0)\Big(A^2+\frac{B^2}{\sin^2 t}+H^2(t)\Big)e^{-2A \int^t_{t_0} \frac{1}{H(\la)} d \la}
\]
\noi or

\ \ \ $\mathrm{(ii)}$  $|f(r)|= r^A$, $|h(\alpha)| =  |\sin\alpha|^A$ and $|g(\theta)| = |g(\theta_0)| e^{\int^{\theta}_{\theta_0} G(t)d t}$, where $G(t)$ satisfies the following
\[
    t + c = 
\begin{cases}
    - \arctan\frac{G(t)}{A} +  \frac{A-1} {\sqrt{A^2-A}}\arctan\frac{G(t)}{\sqrt{A^2-A}}, & \text{if } A^2-A > 0\\
    \frac{1}{G(t)}, & \text{if } A = 0\\
    -\arctan{G(t)}, & \text{if } A = 1\\
    - \arctan\frac{G(t)}{A} +  \frac{A-1}{2 \sqrt{A-A^2}}\ln{\Big|\frac{G(t)-\sqrt{A-A^2}}{G(t)+\sqrt{A-A^2}}\Big|},              & \text{if } A^2-A < 0\\
\end{cases}
\]
\noi or

\ \ \  $\mathrm{(iii)}$  $|g(\theta)| = e^{B\theta}$ , $|f(r)| = |f(r_0)| e^{\int^r_{r_0} \frac{\Phi(t)}{t} d t}$, where $B$ is constant and $\Phi(t)$ satisfies the following
\[
    \ln |t| + c = 
\begin{cases}
   \frac{1}{2}\ln\Big|\frac{\Phi^2(t)+C^2}{\Phi^2(t) - \Phi(t) + C^2}\Big| - \frac{1}{2}\frac{1}{\sqrt{C^2-\frac{1}{4}}}\arctan\frac{\Phi(t)-\frac{1}{2}}{\sqrt{C^2-\frac{1}{4}}}, & \text{if } C^2-\frac{1}{4} > 0\\
   \frac{1}{2}\ln\Big|\frac{\Phi^2(t)+C^2}{\Phi^2(t) - \Phi(t) + C^2}\Big| + \frac{1}{2} \frac{1}{\Phi(t)-\frac{1}{2}}, & \text{if } C^2-\frac{1}{4} = 0\\
   \frac{1}{2}\ln\Big|\frac{\Phi^2(t)+C^2}{\Phi^2(t) - \Phi(t) + C^2}\Big| - \frac{1}{4\sqrt{\frac{1}{4}-C^2}}\ln\Big|\frac{\Phi(t)-\frac{1}{2}-\sqrt{\frac{1}{4} - C^2}}{\Phi(t)-\frac{1}{2}+\sqrt{\frac{1}{4} - C^2}}\Big| & \text{if } C^2-\frac{1}{4} < 0,\\
\end{cases}
\]
and 
\[
|h(\alpha)| = |h(\alpha_0)| e^{B\arcsin\frac{B\cot(\alpha)}{\sqrt{C^2-B^2}}-C\arctan\frac{C\cot\alpha}{\sqrt{C^2-\frac{B^2}{\sin^2 \alpha}}}}
\]
or
\[
|h(\alpha)| = |h(\alpha_0)| e^{-B\arcsin\frac{B\cot(\alpha)}{\sqrt{C^2-B^2}}+C\arctan\frac{C\cot\alpha}{\sqrt{C^2-\frac{B^2}{\sin^2 \alpha}}}}
\]
where c is any constant.
\end{theorem}

\begin{theorem}[Separated two-dimensional $\infty$-Harmonic functions] \label{theorem2} Let $u :\Om \sub\R^2 \larrow \R$ be a $separated$  $\infty$-$harmonic function$ of the $\infty$-Laplace equation
\beq\label{1.14}
u^2_x u_{xx}+2 u_x u_y u_{xy} + u_y^2 u_{yy}=0
\eeq
of the form $u(x,y) = f(x)g(y)$, where $f,g$ are non-linear. Then, one of the following holds: either

\ \ \ $\mathrm{(i)}$  $|f(x)| = |f(x_0)| e^{A (x-x_0)}$ and $|g(y)| = |g(y_0)| e^{\int^y_{y_0}G(t) \,d t}$, where $G(t)$ satisfies 
\[
    t + c = 
\begin{cases}
    \frac{1}{G(t)}, & \text{if } A = 0\\
    -\frac{1}{2 A}\arctan\frac{G(t)}{A}+\frac{G(t)}{2\big(A^2+G^2(t)\big)}, & \text{otherwise } \\
\end{cases}
\]

\noi or

\ \ \ $\mathrm{(ii)}$  $|f(x)| = |f(x_0)| e^{\int^x_{x_0}F(t)\,d t}$  and $|g(y)| = |g(y_0)| e^{B(y-y_0)}$ , where $F(t)$ satisfies

\[
    t + c = 
\begin{cases}
    \frac{1}{F(t)}, & \text{if } B = 0\\
    -\frac{1}{2 B}\arctan\frac{F(t)}{B}+\frac{F(t)}{2\big(B^2+F^2(t)\big)}, & \text{otherwise.} \\
\end{cases}
\]
\end{theorem}

\begin{theorem}[Separated n-dimensional $\infty$-Harmonic functions] \label{theorem3} Let $n\ge 2$ and $u :\Om\sub\R^n \larrow \R$ be a $separated$  $\infty$-$harmonic function$ of the $\infty$-Laplace equation
\beq\label{1.17}
\sum^n_{i,j=1}\D_i u\, \D_j u \, \D^2_{ij}u \, =\, 0.
\eeq
then 
\[
|f_{i}(x_i)|=|f_{i}(x^0_i)|\, e^{A_i(x_i-x^0_i)}\text{ for } 1\le i\neq j \le  n 
\]
and 
\[|f_{j}(x_j)|=|f_{j}(x^0_j)|\, e^{{\displaystyle\int^{x_j}_{x^0_j}F_j(t)\, d t}},
\]
where $F_j(t)$ satisfies
\[
    t + c =  -\frac{1}{2 \Big(\displaystyle\sum_{i\neq j}A^2_i\Big)^{1/2}}\arctan\frac{F_j(t)}{\Big(\displaystyle\sum_{i\neq j}A^2_i\Big)^{1/2}}+\frac{F_j(t)}{2\Big(\displaystyle\sum_{i\neq j}A^2_i+F_j^2(t)\Big)}.
\]
\end{theorem}

\section{Proofs of main results} \label{section2}
In this section we prove our main results. The general idea of our method, which is essentially the same for all our proofs, is to use a substitution to derive a ``better'' PDE. Then, we take any points from the domain which are different only in one component put them to the ``new'' equation and subtract the two equations from each other.

\BPT \ref{theorem1}.
For $u\neq 0$ the equation (\ref{1.1}) can be written as
\beq \label{1.2}
\ \ \frac{u^2_r}{u^2} \frac{u_{rr}}{u} + \frac{2}{r^2}\frac{u_r}{u} \frac{u_\theta}{u} \frac{u_{r\theta}}{u} + \frac{1}{r^4} \frac{u^2_\theta}{u^2} \frac{u_{\theta\theta}}{u} - 
\frac{1}{r^3}\frac{u_r}{u} \frac{u^2_\theta}{u^2} = 0.
\eeq

Let $F=\frac{u_r}{u}$ and $G=\frac{u_\theta}{u}$, 
then $F_r+F^2=\frac{u_{rr}}{u}$, $G_\theta+G^2=\frac{u_{\theta\theta}}{u}$ and $\frac{1}{2}F_\theta+\frac{1}{2}G_r+F G =\frac{u_{r\theta}}{u}$.
Note that $u(r,\theta) = f(r) g(\theta)$ hence $F$ and $G$ doesn't depend on $\theta$ and $r$ respectively, since $F(r,\theta) = \frac{f'(r)}{f(r)}$ and $G(r,\theta) = \frac{g'(\theta)}{g(\theta)}$.
Thus (\ref{1.2}), becomes
\beq \label{1.4}
F^2 F_r+ F^4+ \frac{2}{r^2}F^2 G^2+\frac{1}{r^4}G^2G_{\theta}+\frac{1}{r^4} G^4 - \frac{1}{r^3}F G^2=0.
\eeq
Set $\Phi = F r$, then $r \Phi_r-\Phi = F_r r^2$. Multiplying (\ref{1.4}) by $r^4$, we have
\beq \label{1.6}
(\Phi^2+G^2)(\Phi^2+G^2-\Phi)+r\Phi^2\Phi_r+G^2G_{\theta}=0.
\eeq

We have the following 4 cases for the functions $\Phi$ and $G$:\\
Case (A) $\Phi$ and $G$ are constant functions.\\
Case (B) $\Phi$ is constant and $G$ is non-constant functions.\\
Case (C) $\Phi$ is non-constant and $G$ is constant functions.\\
Case (D) $\Phi$ and $G$ are non-constant functions.\\

Case (A)
Let $\Phi \equiv A$ and $G \equiv  B$, then (\ref{1.6}) gives $A\equiv B\equiv 0$ or $A^2-A+B^2=0$ which can be rewritten as $(A-\frac{1}{2})^2+B^2=\frac{1}{4}$
and as the consequent of substitutions $f(r)=r^A$ and $g(\theta)=e^{B\theta}$ up to a constants.

Case (B)
Let $\Phi \equiv A$, then $G$ is satisfying (\ref{1.6})
\beq\label{1.3}
(A^2+G^2)(A^2+G^2-A)+G^2G_{\theta}=0.
\eeq
Therefore
\[
(A^2+G^2)(A^2+G^2-A)=-G^2 \frac{d G}{d\theta}.\\
\]
Consequently
\beq
\begin{split}
\int{d\theta} &= \int{\frac{-G^2 }{(A^2+G^2)(A^2+G^2-A)}}\,d G\\
&= \int \frac{-A}{A^2+G^2}\,d G - \int \frac{1-A}{A^2 - A+G^2}\,d G.
\end{split}
\eeq
Hence
\[
    t + c = 
\begin{cases}
    - \arctan\frac{G(t)}{A} +  \frac{A-1} {\sqrt{A^2-A}}\arctan\frac{G(t)}{\sqrt{A^2-A}}, & \text{if } A^2-A > 0\\
    \frac{1}{G(t)}, & \text{if } A = 0\\
    -\arctan{G(t)}, & \text{if } A = 1\\
    - \arctan\frac{G(t)}{A} +  \frac{A-1}{2 \sqrt{A-A^2}}\ln{\Big|\frac{G(t)-\sqrt{A-A^2}}{G(t)+\sqrt{A-A^2}}\Big|},              & \text{if } A^2-A < 0.\\
\end{cases}
\]
If A is equal to $0$ or $1$, then $G=0$ is also a solution, hence $g\equiv c$ and $u$ is a linear solution.

Case (C)
Let $G\equiv B$, then $\Phi$ is satisfying (\ref{1.6})
\beq\label{1.7}
(\Phi^2+B^2)(\Phi^2+B^2-\Phi)+r\Phi^2\Phi_r=0.
\eeq
Therefore
\[
(\Phi^2+B^2)(\Phi^2+B^2-\Phi)=-r\Phi^2\frac{d \Phi}{d r}.\\
\]
Consequently
\beq
\begin{split}
\int\frac{1}{r}\,d r &= \int{\frac{-\Phi^2 }{(\Phi^2+B^2)(\Phi^2-\Phi +B^2)}}\,d \Phi \\
&= \int \frac{\Phi}{\Phi^2+B^2}\,d \Phi - \int \frac{\Phi - \frac{1}{2}}{(\Phi-\frac{1}{2})^2 +B^2 -\frac{1}{4}}\,d \Phi -\int \frac{\frac{1}{2}}{(\Phi-\frac{1}{2})^2 +B^2 -\frac{1}{4}}\,d \Phi.\\
\end{split}
\eeq
Hence
\[
    \ln |t| + c = 
\begin{cases}
   \frac{1}{2}\ln\Big|\frac{\Phi^2(t)+B^2}{\Phi^2(t) - \Phi(t) + B^2}\Big| - \frac{1}{2}\frac{1}{\sqrt{B^2-\frac{1}{4}}}\arctan\frac{\Phi(t)-\frac{1}{2}}{\sqrt{B^2-\frac{1}{4}}}, & \text{if } B^2-\frac{1}{4} > 0\\
   \frac{1}{2}\ln\Big|\frac{\Phi^2(t)+B^2}{\Phi^2(t) - \Phi(t) + B^2}\Big| + \frac{1}{2} \frac{1}{\Phi(t)-\frac{1}{2}}, & \text{if } B^2-\frac{1}{4} = 0\\
   \frac{1}{2}\ln\Big|\frac{\Phi^2(t)+B^2}{\Phi^2(t) - \Phi(t) + B^2}\Big| - \frac{1}{4\sqrt{\frac{1}{4}-B^2}}\ln\Big|\frac{\Phi(t)-\frac{1}{2}-\sqrt{\frac{1}{4} - B^2}}{\Phi(t)-\frac{1}{2}+\sqrt{\frac{1}{4} - B^2}}\Big|, & \text{if } B^2-\frac{1}{4} < 0.\\
\end{cases}
\]
If $B=0$, then $\Phi = 0$ is also a solution, hence $u$ is a constant.

Case (D)
Let $\Phi$ and $G$ are non-constant functions, then there exist $r_1\neq r_2$ and  $\theta_1\neq\theta_2$ such that $\Phi(r_1)\neq \Phi(r_2)$ and $G(\theta_1)\neq G(\theta_2)$ satisfying (\ref{1.6}). Thus
\beq \label{1.11}
r_1\Phi(r_1)^2\Phi_r(r_1)-\Phi^3(r_1)+\Phi(r_1)^4+2\Phi(r_1)^2 G(\theta)^2+G(\theta)^2 G_{\theta}(\theta) +G(\theta)^4-\Phi(r_1) G(\theta)^2=0
\eeq
\beq \label{1.12}
r_2\Phi(r_2)^2\Phi_r(r_2)-\Phi^3(r_2)+\Phi(r_2)^4+2\Phi(r_2)^2 G(\theta)^2+G(\theta)^2 G_{\theta}(\theta) +G(\theta)^4-\Phi(r_2) G(\theta)^2=0.
\eeq
Subtracting (\ref{1.11}) and (\ref{1.12}) we get for any $\theta$
\beq
\begin{split}
G^2(\theta) (\Phi(r_1)-\Phi(r_2)) (2(\Phi(r_1)+\Phi(r_2))-1) &= r_2\Phi^2(r_2)\Phi_r(r_2)-\Phi(r_2)^3+\Phi(r_2)^4 \\
& - r_1\Phi^2(r_1)\Phi_r(r_1) + \Phi(r_1)^3 - \Phi(r_1)^4,
\end{split}
\eeq
if $2(\Phi(r_1)+\Phi(r_2))-1\neq 0$, then $G^2(\theta)$ is a constant function or $G(\theta)$ is a step function because $G(\theta_1)=-G(\theta_2)$ and $G(\theta_1)\neq G(\theta_2)$, otherwise $2(\Phi(r_1)+\Phi(r_2))-1 = 0$ for any  $r_1\neq r_2$ such that $\Phi(r_1)\neq \Phi(r_2)$, which means $\Phi(r)$ is a step function and the image of $\Phi$ is symmetrical to $y=\frac{1}{4}$ since $\Phi(r_1)-\frac{1}{4} = -(\Phi(r_2)-\frac{1}{4})$. For both cases we have a contradiction to $C^{1,\alpha}$ regularity for $\infty$-Harmonic mappings in two dimensions (\cite{ES}, \cite{S}).

Finally integrating and substituting we complete the proof.
\qed

\begin{remark}[The Arronson solution]
\label{Remark 1} 
If $A=\frac{4}{3}$ then $A^2-A>0$ and $G(t)$ function satisfies
\[
\begin{split}
t +c &= - \arctan\frac{3}{4}G(t) +  \frac{1} {2}\arctan\frac{3}{2}G(t),\\
\end{split}
\]
which can be rewritten as
\[
27 G^3(t)+54G^2(t)\tan2(t+c)+32\tan2(t+c)=0.
\]
Solving a third degree equation with respect to $G(t)$, we get
\[
G(t)=-\frac{4}{3}\frac{\tan^{\frac{1}{3}}(t+c)+\tan^{\frac{5}{3}}(t+c)+\tan(t+c)}{1-\tan^2(t+c)}.
\]
Therefore
\[
\begin{split}
\int G(t)\,d t&=
 \ln\Bigg|\frac{\big(1-\tan^{\frac{2}{3}}(t+c)\big)\big(1+\tan^\frac{2}{3}(t+c)\big)^\frac{1}{3}}{\big(\tan^\frac{4}{3}(t+c)-\tan^\frac{2}{3}(t+c)+1\big)^\frac{2}{3}}\Bigg|.
\end{split}
\]
Hence
\[
\begin{split}
e^{\int^{\theta}_{\theta_0} G(t)\, d t} &=\frac{|1-\tan^\frac{2}{3}(\theta+c)| |1+\tan^\frac{2}{3}(\theta+c)|^\frac{1}{3}}{|\tan^\frac{4}{3}(\theta+c)-\tan^\frac{2}{3}(\theta+c)+1|^\frac{2}{3}}  \cdot c(\theta_0)\\
& = \frac{|1-\tan^\frac{4}{3}(\theta+c)| }{|1+\tan^2(\theta+c)|^\frac{2}{3}}\cdot c(\theta_0)\\
& = \big|\cos^\frac{4}{3}(\theta+c)-\sin^\frac{4}{3}(\theta+c)\big|  \cdot c(\theta_0).
\end{split}
\]
Finally
\[
\begin{split}
|g(\theta)| & = |g(\theta_0)| e^{\int^{\theta}_{\theta_0} G(t)}d t  \\
& = |g(\theta_0)|\big(\big|\cos^\frac{4}{3}(\theta+c)-\sin^\frac{4}{3}(\theta+c)\big| \big)\\
|f(r)| &= r^\frac{4}{3}.\\
\end{split}
\]
Thus, one of the possible solutions
\[
\begin{split}
u(r,\theta) &= f(r) g(\theta)\\
& = r^\frac{4}{3}\big(\cos^\frac{4}{3}(\theta+c)-\sin^\frac{4}{3}(\theta+c)\big).\\
u(x,y) & = |x|^\frac{4}{3}-|y|^\frac{4}{3}.
\end{split}
\]
\end{remark}

\begin{remark}[The Aronsson solution] 
If $A=-\frac{1}{3}$ then $A^2-A>0$, similarly as in previous case we can find that $u(r,\theta) = r^{-\frac{1}{3}}\big(\cos^\frac{4}{3}(\frac{\theta+c}{2})-\sin^\frac{4}{3}(\frac{\theta+c}{2})\big)$ is the solution of the $\infty$-Laplace equation which was described in \cite{A3} as
\[
g(\theta) = \frac{\cos t}{(1+3 \cos^2 t)^\frac{2}{3}}, \ \ \ \ \ \theta=t-2 \arctan\bigg(\frac{\tan t}{2}\bigg), \ \ \ \ \ -\frac{\pi}{2}\le t \le \frac{\pi}{2}.
\]
The key fact two solutions are identically equal  is $\tan\frac{\theta}{2}=-\tan^3\frac{t}{2}$.
\end{remark}

\BPT \ref{theorem2}. It is particular case of the Theorem \ref{theorem3}, when $n = 2$.\qed

\BPT \ref{theorem3}. For $u\neq 0$ equation (\ref{1.17}) can be written as
\beq \label{1.18}
\sum^n_{i,j=1}\frac{\D _i u}{u}\frac{\D _j u}{u} \frac{\D^2_{ij}u}{u}=0.
\eeq
Let $F_i=\frac{\D _i u}{u}$ then $\D_i F_i + F^2_i = \frac{\D_{ii}u}{u}$ and $F_i F_j = \frac{\D_{ij}u}{u}$.
Thus (\ref{1.18}), becomes
\beq \label{1.19}
\Bigg(\sum^n_{i = 1} F_i (x_i) ^2\Bigg)^2+ \sum^n_{i=1} F_i (x_i) ^2 D_i F_i(x_i)=0.
\eeq
Since $u(x)=\displaystyle\prod^n_{i=1} f_i (x_i)$, then $F_i$ depends only on $x_i$, consequently $\D_i F_i(x_i) = F'_i (x_i)$. Set $x^1$, $x^2\in\Om$ such that $x^{1} = (x_1,x_2,...,x^{1}_j,...,x_n)$ and $x^{2} = (x_1,x_2,...,x^{2}_j,...,x_n)$, where $x^1_j\neq x^2_j$ in (\ref{1.19}) and subtract two equations. We find
\[
\Big(F^2_j(x^1_j) - F^2_j(x^2_j)\Big)\Big(2\sum^n_{i\neq j}F^2_i(x_i)+2 F^2_j(x^1_j) + 2 F^2_j(x^1_j)\Big) + F^2_j F'_j(x^1_j) - F^2_j F'_j(x^2_j) = 0,
\]
assuming $F^2_j(x^1_j) \neq F^2_j(x^2_j)$, we have
\beq \label{1.20}
2\sum^n_{i\neq j}F^2_i(x_i) = -\frac{ F^2_j F'_j(x^1_j) - F^2_j F'_j(x^2_j)}{F^2_j(x^1_j) - F^2_j(x^2_j)}-2 F^2_j(x^1_j) - 2 F^2_j(x^1_j).
\eeq
LHS of ($\ref{1.20}$) does not depend on $x^1_j$ and $x^2_j$ so
\[
\sum^n_{i\neq j}F^2_i(x_i) \equiv c
\]
for all $x_i$. Hence $F_i(x_i)=A_i$, where $A_i$ is a constant for all $i\neq j$. Thus ($\ref{1.19}$) gives
\[
\Bigg(\sum^n_{i = 1}A^2_i+ F_j (x_j) ^2\Bigg)^2+ F_j (x_j) ^2 F'_j(x_j)=0,
\]
consequently
\[
d x_j = -\frac{F^2_j}{\Bigg(\displaystyle\sum^n_{i \neq j}A^2_i+ F_j^2\Bigg)^2} d F_j,
\]
hence
\[
x_j + c =  -\frac{1}{2 \sqrt{\displaystyle \sum_{i\neq j}A^2_i}}\arctan\frac{F_j(x_j)}{\sqrt{\displaystyle\sum_{i\neq j}A^2_i}}+\frac{F_j(x_j)}{2\Big(\displaystyle\sum_{i\neq j}A^2_i+ F_j^2(x_j)\Big)}, \text{ if } \displaystyle\sum^n_{i\neq j}A^2_i \neq 0.
\]
Otherwise 
\[
F_j^4 (x_j)+F^2(x_j)F'(x_j) = 0,
\]
so
\beq\label{1.21}
F^2(x_j)+F'_j(x_j) = 0, \text{ since we assume } F^2_j(x^1_j) \neq F^2_j(x^2_j).
\eeq
Solving (\ref{1.21})  we get $F_j(x_j) = -\frac{1}{x_j + c}$. Hence $f_i(x_i) = x_i + c_i$ for all $i \neq j$ and $f_j(x_j) = \frac{c_j}{|x_j+c|}$, where $c_i$ are constants for all $1 \leq i \leq n$.

If there is no $j$ such that $F^2_j(x^1_j) \neq F^2_j(x^2_j)$ then $F^2_j(x_j) \equiv c_j$ for all $1\leq j \leq n$ and (\ref{1.19}) gives that all $c_j=0$ for all $1\leq j\leq n$. So $f_i(x_i)=C_i$, where $C_i$ is constant for all $i$.
\qed

\BPT \ref{theorem4}. For $u\neq 0$ we can rewrite (\ref{1.16}) as 
\beq
\begin{split}
& \frac{u^2_r}{u^2} \frac{u_{r r}}{u} + \frac{2}{r^2 \sin^2\alpha}\frac{u_r}{u} \frac{u_\theta}{u}  \frac{u_{r\theta}}{u} + \frac{1}{r^4\sin^4\alpha} \frac{u^2_\theta}{u^2} \frac{u_{\theta\theta}}{u} - \frac{1}{r^3\sin^2\alpha}\frac{u_r}{u} \frac{u^2_\theta}{u^2} + \frac{2}{r^2}\frac{u_r}{u} \frac{u_{\alpha}}{u} \frac{u_{r \alpha}}{u}+\\
\ \ \ &  + \frac{1}{r^4}\frac{u^2_{\alpha}}{u^2} \frac{u_{\alpha\alpha}}{u} - \frac{1}{r^3}\frac{u_r}{u} \frac{u^2_{\alpha}}{u^2}+\frac{2}{r^4 \sin^2 \alpha}\frac{u_{\theta}}{u} \frac{u_{\alpha}}{u} \frac{u_{\theta\alpha}}{u}-\frac{\cos\alpha}{r^4\sin^3\alpha} \frac{u_{\alpha}}{u}\frac{u^2_\theta}{u^2}= 0.\\
\end{split}
\eeq
Let 
\[
\text{$F = \frac{u_r}{u}$, \ $G = \frac{u_\theta}{u}$ \ and \ $H = \frac{u_\alpha}{u}$}
\]
Then, we have 
\[
\text{$F' + F^2 = \frac{u_{r r}}{u}$, \ $G' + G^2 = \frac{u_{\theta\theta}}{u}$ \ and \ $H' + H^2 = \frac{u_{\alpha\alpha}}{u}$.}
\]
Also, 
\[
\text{$\frac{u_{r\theta}}{u} = F' G'$, \ $\frac{u_{r\alpha}}{u} = F' H'$ \ and \ $\frac{u_{\theta\alpha}}{u} = G' H'$.}
\] 
Note that $F$, $G$ and $H$ depend on only $r$, $\theta$ and $\alpha$ respectively, since $\frac{u_r}{u} = \frac{f'}{f}$, $\frac{u_\theta}{u} = \frac{g'}{g}$ and $\frac{u_\alpha}{u} = \frac{h'}{h}$. Thus (\ref{1.21}) becomes 
\[
\begin{split}
&F^2(F'+F^2)+ \frac{2}{r^2 \sin^2\alpha}F^2G^2+ \frac{1}{r^4\sin^4\alpha}G^2(G^2+G')- \frac{1}{r^3\sin^2\alpha}F G^2+\\
&+ \frac{2}{r^2}F^2 H^2+\frac{1}{r^4}H^2(H^2 + H')-\frac{1}{r^3}F H^2 + \frac{2}{r^4\sin^2\alpha}G^2H^2 - \frac{\cos\alpha}{r^4\sin^3\alpha}HG^2=0,\\
\end{split}
\]
which is equivalent to
\beq\label{1.22}
\begin{split}
&\Big(F^2+\frac{1}{r^2\sin^2\alpha}G^2+\frac{1}{r^2}H^2\Big)^2+F^2F'+\frac{1}{r^4\sin^4\alpha}G^2G'+\frac{1}{r^4}H^2H'-\\
\ \ \ \ &- \frac{\cos\alpha}{r^4\sin^3\alpha}HG^2 - \frac{1}{r^3\sin^2\alpha}F G^2-\frac{1}{r^3}F H^2=0.\\
\end{split}
\eeq
Let $\Phi = r F$. Then, $r\Phi'-\Phi = r^2F'$. Multiplying (\ref{1.22}) by $r^4$ we have
\beq\label{1.23}
\begin{split}
&\Big(\Phi(r)^2+\frac{1}{\sin^2\alpha}G(\theta)^2+H(\alpha)^2\Big) \Big(\Phi(r)^2-\Phi(r)+\frac{1}{\sin^2\alpha}G(\theta)^2 + H(\alpha)^2 \Big)- \\
& \frac{\cos\alpha}{\sin^3\alpha}H(\alpha)G(\theta)^2 + r\Phi(r)^2\Phi'(r) + \frac{1}{\sin^4\alpha}G(\theta)^2 G'(\theta) + H(\alpha)^2 H'(\alpha)=0.\\
\end{split}
\eeq

Setting $(r_1,\theta,\alpha)$, $(r_2,\theta,\alpha)\in\Om$ such that $r_1\neq r_2$ in (\ref{1.23}) and subtracting two equation we get
\[
\begin{split}
&\Big(2\Phi^2(r_1) - \Phi(r_1) - 2\Phi^2(r_2) + \Phi(r_2) \Big) \Big(H^2(\alpha)+ \frac{1}{\sin^2\alpha} G^2(\theta)\Big)+\\
&\Phi^4(r_1)-\Phi^4(r_2)-\Phi^3(r_1)+\Phi^3(r_2)+ r_1\Phi^2(r_1)\Phi'(r_1) - r_2\Phi^2(r_2)\Phi'(r_2) = 0.\\
\end{split}
\]

Case 1. If $2\Phi^2(r)-\Phi(r)\not\equiv c$, then $H^2(\alpha)+ \frac{1}{\sin^2\alpha} G^2(\theta)= C^2$, hence $G(\theta)\equiv B$, $H(\alpha)=\pm\sqrt{C^2-\frac{B^2}{\sin^2\alpha}}$
and (\ref{1.23}) becomes
\beq\label{1.24}
\Big(\Phi^2(r) + C^2\Big)\Big(\Phi^2(r) - \Phi(r) + C^2\Big)+r\Phi^2(r)\Phi'(r)= H(\alpha)B^2\frac{\cos\alpha}{\sin^3\alpha}-H^2(\alpha)H'(\alpha).
\eeq
LHS in (\ref{1.24}) depends on $r$ only and RHS in (\ref{1.24}) depends on $\alpha$ only, so we have the system
\beq\label{1.25}
\begin{cases}
\Big(\Phi^2(r) + C^2\Big)\Big(\Phi^2(r) - \Phi(r) + C^2\Big)+r\Phi^2(r)\Phi'(r) =c_1 \\
H(\alpha)B^2\frac{\cos\alpha}{\sin^3\alpha}-H^2(\alpha)H'(\alpha) = c_1\\
H(\alpha)=\pm\sqrt{C^2-\frac{B^2}{\sin^2\alpha}}\\
G(\theta)\equiv B.
\end{cases}
\eeq
Solving (\ref{1.25}) we get $c_1=0$, $G(\theta)\equiv B$, $H(\alpha)=\pm\sqrt{C^2-\frac{B^2}{\sin^2\alpha}}$ and from the previous result (\ref{1.7}) $\Phi$ satisfies the following
\[
    \ln |t| + c = 
\begin{cases}
   \frac{1}{2}\ln\Big|\frac{\Phi^2(t)+C^2}{\Phi^2(t) - \Phi(t) + C^2}\Big| - \frac{1}{2}\frac{1}{\sqrt{C^2-\frac{1}{4}}}\arctan\frac{\Phi(t)-\frac{1}{2}}{\sqrt{C^2-\frac{1}{4}}}, & \text{if } C^2-\frac{1}{4} > 0\\
   \frac{1}{2}\ln\Big|\frac{\Phi^2(t)+C^2}{\Phi^2(t) - \Phi(t) + C^2}\Big| + \frac{1}{2} \frac{1}{\Phi(t)-\frac{1}{2}}, & \text{if } C^2-\frac{1}{4} = 0\\
   \frac{1}{2}\ln\Big|\frac{\Phi^2(t)+C^2}{\Phi^2(t) - \Phi(t) + C^2}\Big| - \frac{1}{4\sqrt{\frac{1}{4}-C^2}}\ln\Big|\frac{\Phi(t)-\frac{1}{2}-\sqrt{\frac{1}{4} - C^2}}{\Phi(t)-\frac{1}{2}+\sqrt{\frac{1}{4} - C^2}}\Big|, & \text{if } C^2-\frac{1}{4} < 0.\\
\end{cases}
\]
If $B = 0$, then $C = 0$, $\Phi=0$ is also a solution of (\ref{1.25}).

Case 2. If $2\Phi^2(r)-\Phi(r)\equiv c$, then $\Phi(r)\equiv A$. Setting $(r,\theta_1,\alpha)$, $(r,\theta_2,\alpha)\in\Om$ such that $\theta_1\neq \theta_2$ in (\ref{1.23}) and subtracting two equations, we get
\[
\begin{split}
&\Big(G^2(\theta_1)-G^2(\theta_2)\Big)\Big(2H^2(\alpha)-\cot(\alpha) H(\alpha)+2A^2-A\Big)\frac{1}{\sin^2\alpha}+\\
&+\frac{1}{\sin^4\alpha}\Big(G^2(\theta_1)G'(\theta_1)-G^2(\theta_2)G(\theta_2)+G^4(\theta_1)-G^4(\theta_2)\Big) = 0.
\end{split}
\]

Case 2a. If $G^2(\theta)\not\equiv c$, then $\Big(2H^2(\alpha)-\cot(\alpha) H(\alpha)+2A^2-A\Big)\sin^2\alpha \equiv C$ and (\ref{1.23}) becomes
\[
\begin{split}
&\frac{G^4(\theta)}{\sin^4\alpha}+\frac{G^2(\theta)}{\sin^2\alpha}\Big(2H^2(\alpha)-\cot(\alpha)H(\alpha)+2\Phi(r)^2-\Phi(r)\Big)+ \frac{G^2(\theta)G'(\theta)}{\sin^4\alpha} = \\
& \Phi^3(r)-\Phi^4(r)-2H^2(\alpha)\Phi^2(\theta)+\Phi(r)H^2(\alpha)-H^4(\alpha)-H^2(\alpha)H'(\alpha),\\
\end{split}
\]
thus
\beq\label{1.26}
\begin{split}
&G^4(\theta)+G^2(\theta)C+ G^2(\theta)G'(\theta) = \sin^4\alpha\Big(\Phi^3(r)-\Phi^4(r)-2H^2(\alpha)\Phi^2(\theta)+\\
&\Phi(r)H^2(\alpha)-H^4(\alpha)-H^2(\alpha)H'(\alpha)\Big).
\end{split}
\eeq

The LHS of (\ref{1.26}) depends on $\theta$ only and the RHS depends on $r$ and $\alpha$ only, hence we have a system
\beq \label{1.27}
\begin{cases}
G^4(\theta)+G^2(\theta) C+ G^2(\theta)G'(\theta) = c_1 \\
\sin^4\alpha\Big(\Phi^3(r)-\Phi^4(r)-2H^2(\alpha)\Phi^2(\theta)+\Phi(r)H^2(\alpha)-H^4(\alpha)-H^2(\alpha)H'(\alpha)\Big) = c_1\\
\Phi(r)\equiv A\\
\Big(2H^2(\alpha)-\cot(\alpha) H(\alpha)+2A^2-A\Big)\sin^2\alpha \equiv C.\\
\end{cases}
\eeq

Solving (\ref{1.27}) we  get $C = 2A^2-A$, $c_1 = - A^4 + A^3$, $\Phi(r) = A$, $H(\alpha) = A\cot(\alpha)$ and $G(\theta)$ satisfies the following
\[
\Big(G^2(\theta)+A^2\Big)\Big(G^2(\theta)+A^2-A\Big) + G^2(\theta)G'(\theta) = 0,
\]
which is equivalent from the previous result (\ref{1.3}) to
\[
    t + c = 
\begin{cases}
    - \arctan\frac{G(t)}{A} +  \frac{A-1} {\sqrt{A^2-A}}\arctan\frac{G(t)}{\sqrt{A^2-A}}, & \text{if } A^2-A > 0\\
    \frac{1}{G(t)}, & \text{if } A = 0\\
    -\arctan{G(t)}, & \text{if } A = 1\\
    - \arctan\frac{G(t)}{A} +  \frac{A-1}{2 \sqrt{A-A^2}}\ln{\Big|\frac{G(t)-\sqrt{A-A^2}}{G(t)+\sqrt{A-A^2}}\Big|},              & \text{if } A^2-A < 0.\\
\end{cases}
\]
If A is equal to $0$ or $1$, then $G=0$ is also a solution of (\ref{1.27}).

Case 2b. If $G^2(\theta)\equiv c$, then $G(\theta) = B$. Thus (\ref{1.23}) gives 
\beq\label{1.28}
\Big(A^2-A+\frac{B^2}{\sin^2\alpha}+H^2(\alpha)\Big)\Big(A^2+\frac{B^2}{\sin^2\alpha}+H^2(\alpha)\Big)-\frac{B^2\cos\alpha}{\sin^3\alpha}H(\alpha)+H^2(\alpha)H'(\alpha) = 0.
\eeq
Let $\Psi(\alpha) = \frac{B^2}{\sin^2\alpha}+H^2(\alpha)$, then (\ref{1.28}) becomes 
\[
(A^2-A+\Psi)(A^2+\Psi)+ H(\alpha)\frac{1}{2}\Psi'(\alpha) = 0,
\]
which gives us
\[
1+\frac{H}{2A}\Bigg(\ln\Bigg|\frac{A^2-A+\Psi}{A^2+\Psi}\Bigg|\Bigg)_\alpha = 0,
\]
consequently
\[
|A^2-A+\Psi(t)| = c(A,\Psi(t_0))(A^2+\Psi(t)) e^{-2A\int^t_{t_0}\frac{1}{H(\la)} d \la},
\]
so that
\[
\Bigg|A^2-A+\frac{B^2}{\sin^2t}+H^2(t)\Bigg| = c(A, B,H(t_0))\Bigg(A^2+\frac{B^2}{\sin^2t}+H^2(t)\Bigg) e^{\displaystyle{-2A\int^t_{t_0}\frac{1}{H(\la)} d \la}}.
\]

Finally integrating and substituting we complete the proof.
\qed

\section{Numerical approximations of $\infty$-harmonic functions} \label{section3}

In this section we illustrate the $\infty$-Harmonic functions derived earlier, depending on the parameter(s). The results illustrate that we may have a family of solutions depending on the 2$\pi$-interval even if the parameter(s) is/are fixed. For example: the solution on Figure \ref{fig:14} is a combination of those in Figure \ref{fig:15} and Figure \ref{fig:16} when $\theta$ belongs to 1st and 2nd $2\pi$- interval of the domain respectively.

\begin{figure}[h]
\caption{The approximation to $u$ of the Theorem \ref{theorem1} i, depending on the parameters $A$ and $B$.}\label{fig:thm1typei}
    \centering
    \begin{subfigure}[b]{0.3\textwidth}
        \includegraphics[width=\textwidth]{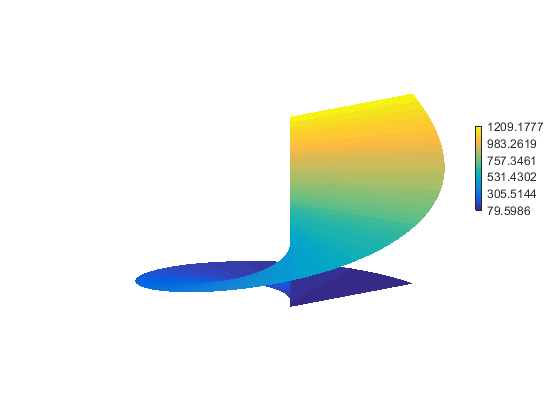}
        \caption{$A = 0.25, B = 0.433$}
        \label{fig:1}
    \end{subfigure}
    \begin{subfigure}[b]{0.3\textwidth}
        \includegraphics[width=\textwidth]{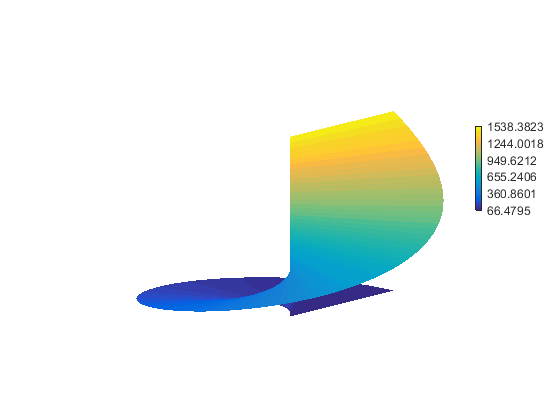}
        \caption{$A = 0.5, B = 0.5$}
        \label{fig:2}
    \end{subfigure}
    \begin{subfigure}[b]{0.3\textwidth}
        \includegraphics[width=\textwidth]{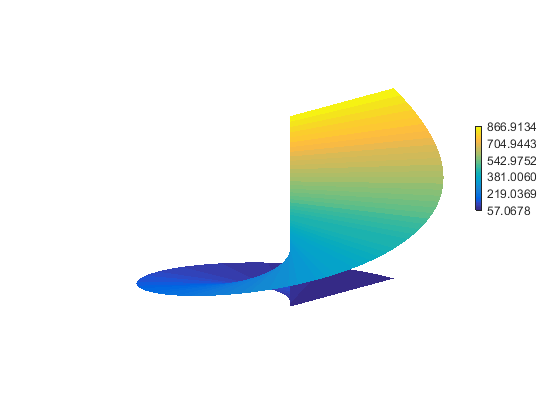}
        \caption{$A = 0.75, B = 0.433$}
        \label{fig:3}
    \end{subfigure}
    
    \begin{subfigure}[b]{0.3\textwidth}
        \includegraphics[width=\textwidth]{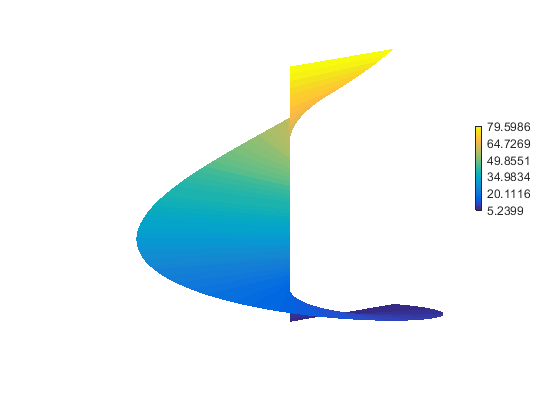}
        \caption{$A = 0.25, B = -0.433$}
        \label{fig:4}
    \end{subfigure}
    \begin{subfigure}[b]{0.3\textwidth}
        \includegraphics[width=\textwidth]{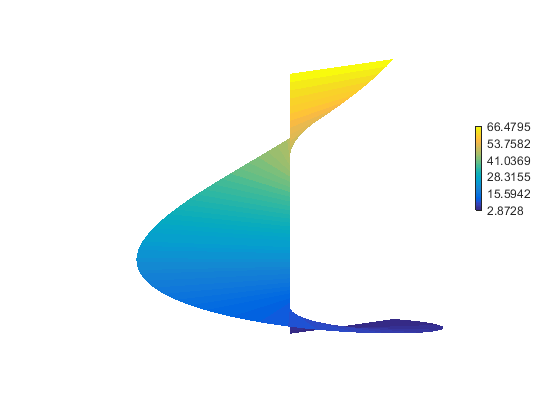}
        \caption{$A = 0.5, B = -0.5$}
        \label{fig:5}
    \end{subfigure}
    \begin{subfigure}[b]{0.3\textwidth}
        \includegraphics[width=\textwidth]{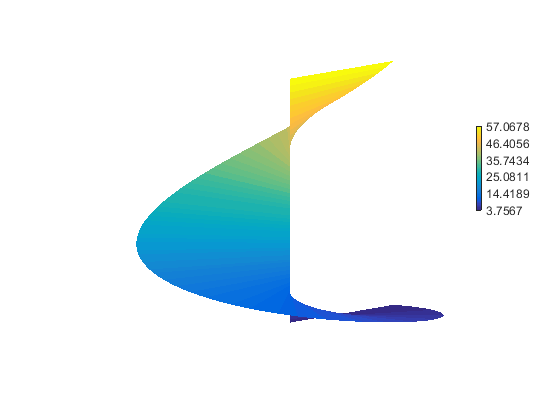}
        \caption{$A = 0.75, B = -0.433$}
        \label{fig:6}
    \end{subfigure}
\end{figure}
\begin{figure}[h]
\captionsetup{type=figure}\addtocounter{figure}{-1}
\caption{The approximation to $u$ of the Theorem \ref{theorem1} ii, depending on the parameter $A$.}\label{fig:thm1typeii}
    \centering
    \begin{subfigure}[b]{0.3\textwidth}
        \includegraphics[width=\textwidth]{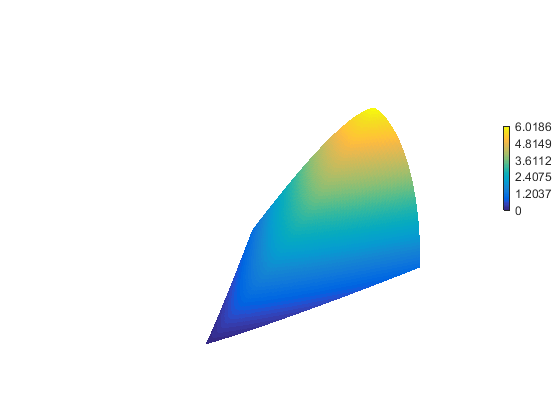}
        \caption{$A = 4/3$}
        \label{fig:7}
    \end{subfigure}
    \begin{subfigure}[b]{0.3\textwidth}
        \includegraphics[width=\textwidth]{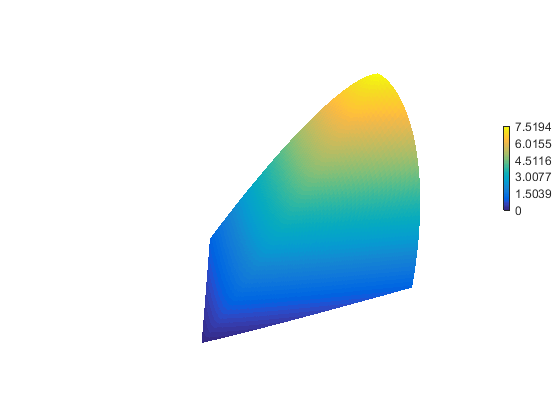}
        \caption{$A = 1.15$}
        \label{fig:8}
    \end{subfigure}
    \begin{subfigure}[b]{0.3\textwidth}
        \includegraphics[width=\textwidth]{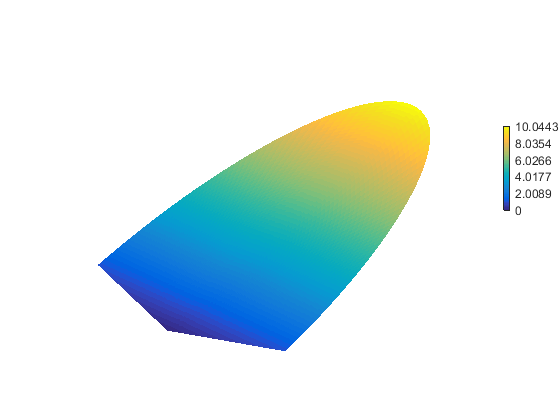}
        \caption{$A = 1$}
        \label{fig:9}
    \end{subfigure}
    
    \begin{subfigure}[b]{0.3\textwidth}
        \includegraphics[width=\textwidth]{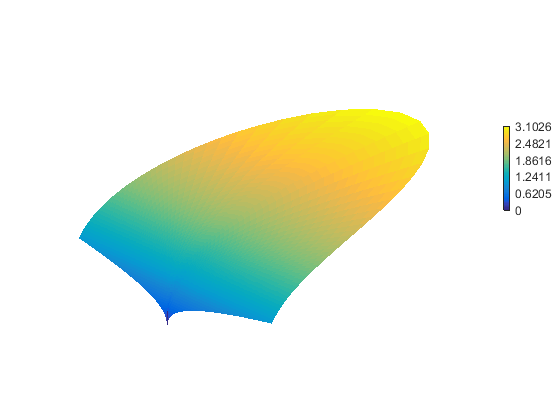}
        \caption{$A = 1/3$}
        \label{fig:10}
    \end{subfigure}
    \begin{subfigure}[b]{0.3\textwidth}
        \includegraphics[width=\textwidth]{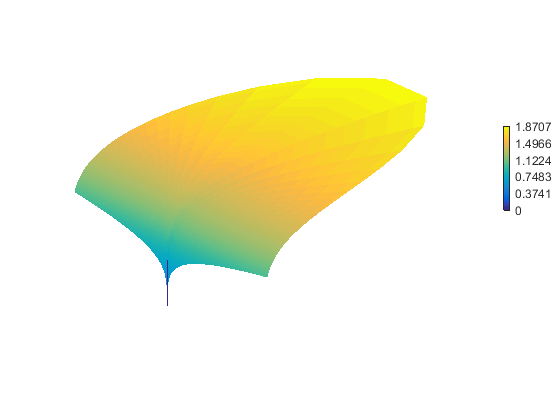}
        \caption{$A = 0.15$}
        \label{fig:11}
    \end{subfigure}
    \begin{subfigure}[b]{0.3\textwidth}
        \includegraphics[width=\textwidth]{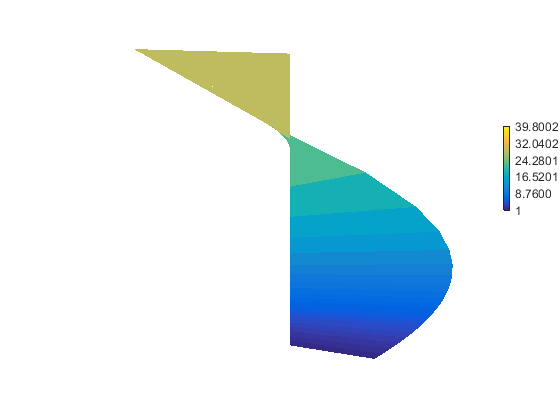}
        \caption{$A = 0$}
        \label{fig:12}
    \end{subfigure}
    
    \begin{subfigure}[b]{0.3\textwidth}
        \includegraphics[width=\textwidth]{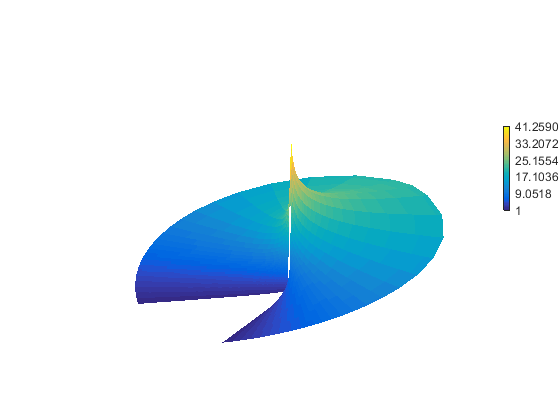}
        \caption{$A = -0.15$}
        \label{fig:13}
    \end{subfigure}
    \begin{subfigure}[b]{0.3\textwidth}
        \includegraphics[width=\textwidth]{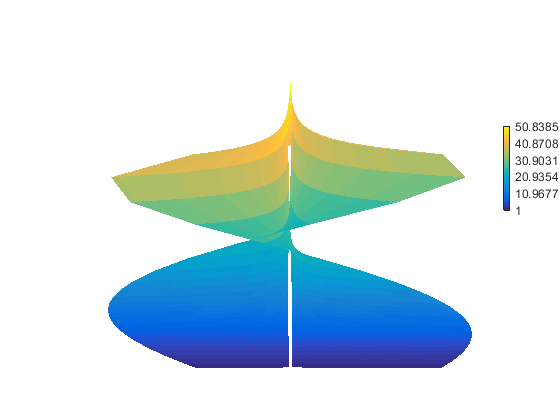}
        \caption{$A =-0.05$}
        \label{fig:14}
    \end{subfigure}

    \begin{subfigure}[b]{0.3\textwidth}
        \includegraphics[width=\textwidth]{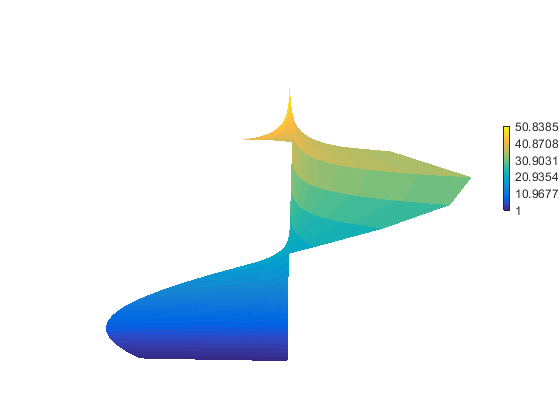}
        \caption{$A =-0.05$}
        \label{fig:15}
    \end{subfigure}
    \begin{subfigure}[b]{0.3\textwidth}
        \includegraphics[width=\textwidth]{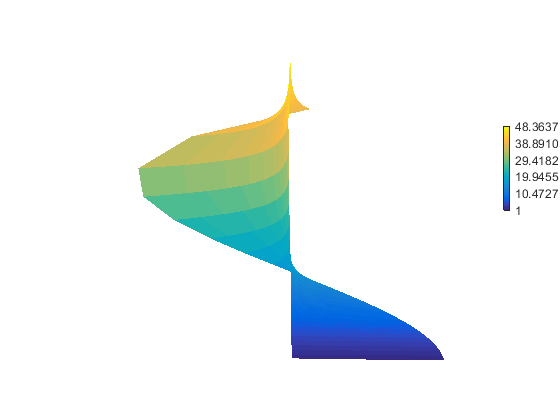}
        \caption{$A = -0.05$}
        \label{fig:16}
    \end{subfigure}
\end{figure}
\begin{figure}[h]
\captionsetup{type=figure}\addtocounter{figure}{-1}
\caption{The approximation to $u$ of the Theorem \ref{theorem1} iii, depending on the parameter $B$.}\label{fig:thm1typeiii}
    \centering
    \begin{subfigure}[b]{0.3\textwidth}
        \includegraphics[width=\textwidth]{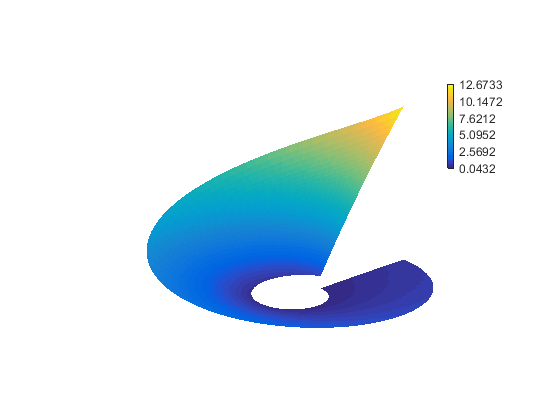}
        \caption{$B = -1/3$}
        \label{fig:17}
    \end{subfigure}
    \begin{subfigure}[b]{0.3\textwidth}
        \includegraphics[width=\textwidth]{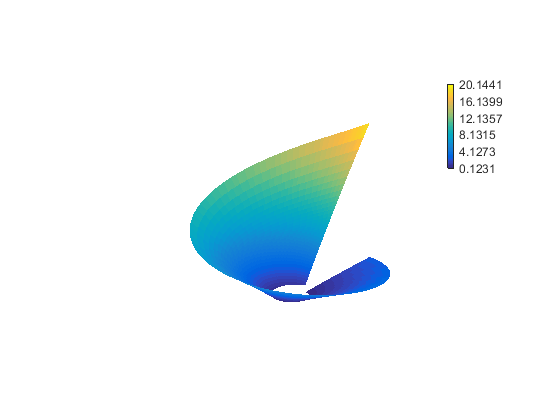}
        \caption{$ B = -1/2$}
        \label{fig:18}
    \end{subfigure}
    \begin{subfigure}[b]{0.3\textwidth}
        \includegraphics[width=\textwidth]{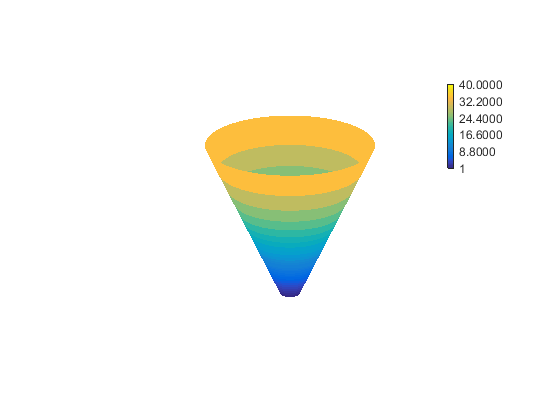}
        \caption{$B = 0$}
        \label{fig:19}
    \end{subfigure}
    
    \begin{subfigure}[b]{0.3\textwidth}
        \includegraphics[width=\textwidth]{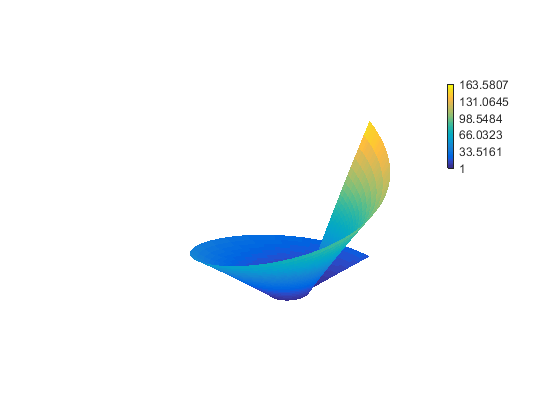}
        \caption{$B = 1/3$}
        \label{fig:20}
    \end{subfigure}
    \begin{subfigure}[b]{0.3\textwidth}
        \includegraphics[width=\textwidth]{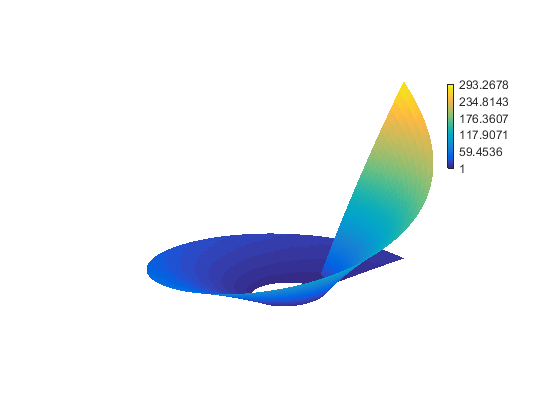}
        \caption{$B =1/2$}
        \label{fig:21}
    \end{subfigure}
    \begin{subfigure}[b]{0.3\textwidth}
        \includegraphics[width=\textwidth]{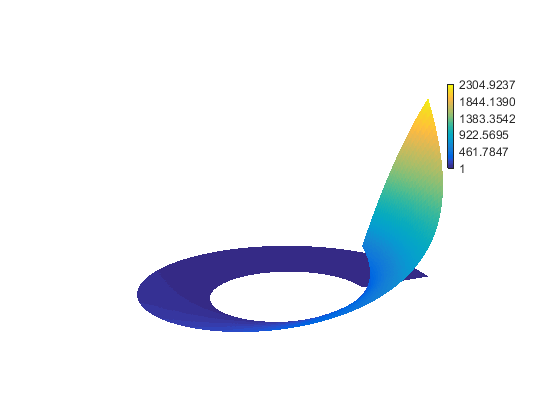}
        \caption{$B = 1$}
        \label{fig:22}
    \end{subfigure}
\end{figure}
\begin{figure}[h]
\captionsetup{type=figure}\addtocounter{figure}{-1}
\caption{The approximation to $u$ of the Theorem \ref{theorem2} i, depending on the parameter $A$.}\label{fig:thm2}
    \centering
    \begin{subfigure}[b]{0.3\textwidth}
        \includegraphics[width=\textwidth]{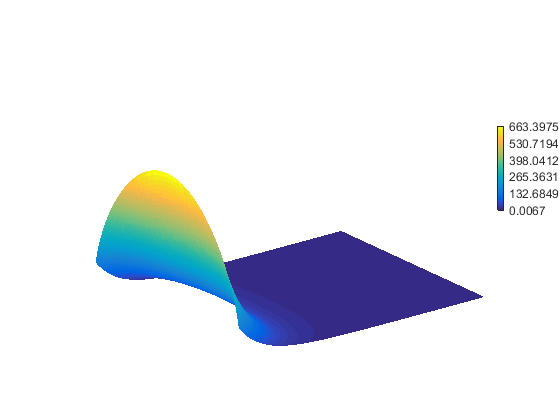}
        \caption{$A = -0.5$}
        \label{fig:23}
    \end{subfigure}
    \begin{subfigure}[b]{0.3\textwidth}
        \includegraphics[width=\textwidth]{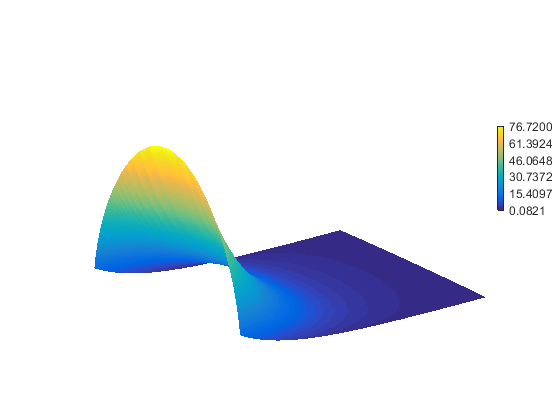}
        \caption{$ A = -0.25$}
        \label{fig:24}
    \end{subfigure}
    \begin{subfigure}[b]{0.3\textwidth}
        \includegraphics[width=\textwidth]{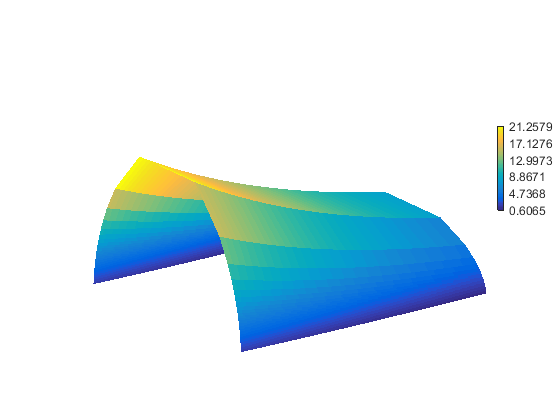}
        \caption{$A = -0.05$}
        \label{fig:25}
    \end{subfigure}

    \begin{subfigure}[b]{0.3\textwidth}
        \includegraphics[width=\textwidth]{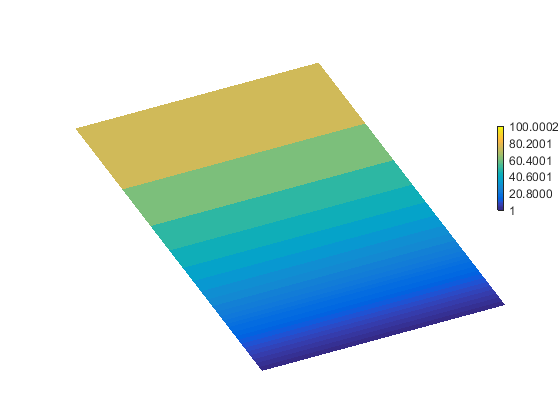}
        \caption{$A = 0$}
        \label{fig:26}
    \end{subfigure}
    \begin{subfigure}[b]{0.3\textwidth}
        \includegraphics[width=\textwidth]{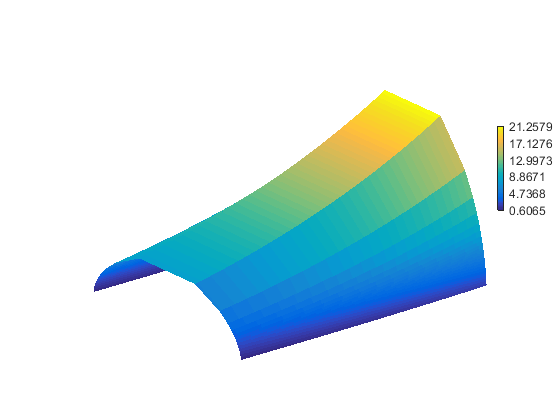}
        \caption{$A = 0.05$}
        \label{fig:27}
    \end{subfigure}
    \begin{subfigure}[b]{0.3\textwidth}
        \includegraphics[width=\textwidth]{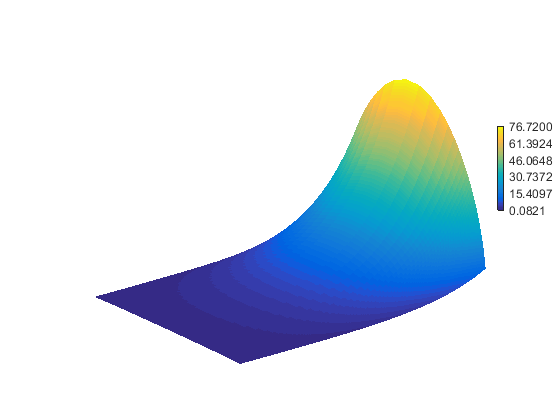}
        \caption{$A = 0.25$}
        \label{fig:28}
    \end{subfigure}
\end{figure}

\ms

\noi \textbf{Acknowledgement.} B.A. would like to thank Nikos Katzourakis for improvements of the presentation of this paper and Tristan Pryer for comments on numerical experiments.

\ms

\bibliographystyle{amsplain}

\begin{thebibliography}{30}

\bibitem[AK]{AK} H. Abugirda, N. Katzourakis, \emph{Existence of 1D Vectorial Absolute Minimisers in $L^\infty$ under Minimal Assumptions}, Proceeding of the AMS 145, 2567 - 2575 (2017).

\bibitem[AyK]{AyK} B. Ayanbayev, N. Katzourakis, \emph{A Pointwise Characterisation of the PDE system of vectorial Calculus of variations in $L^\infty$}, Proc. Royal Soc. Edinburgh A, in press.

\bibitem[A1]{A1} G. Aronsson, \emph{Extension of functions satisfying Lipschitz conditions}, Arkiv f\"ur Mat. 6 (1967), 551 - 561.

\bibitem[A2]{A2} G. Aronsson, \emph{On the partial differential equation $u_x^2 u_{xx} + 2u_x u_y u_{xy} + u_y^2 u_{yy} = 0$}, Arkiv f\"ur Mat. 7 (1968), 395 - 425.

\bibitem[A3]{A3} G. Aronsson, \emph{Construction of singular solutions to the p-harmonic equation and its limit equation for $p=\infty$}, Manuscripta mathematica Volume 56 issue 2 1986, 135 - 158.

\bibitem[CKP]{CKP} G. Croce, N. Katzourakis, G. Pisante, \emph{$\mD$-solutions to the system of vectorial Calculus of Variations in $L^\infty$ via the Baire Category method for the singular values}, ArXiv preprint, https://arxiv.org/pdf/1509.01811.pdf.

\bibitem[ES]{ES} L.C.Evans, O.Savin, \emph{$C^{1,\alpha}$ regularity for $\infty$-Harmonic functions in two dimensions}, Calculus of Variations (2008) 32, 325-347.

\bibitem[FF]{FF} I.L. Freire, A. C. Faleiros, \emph{Lie point symmetries and some group invariant solutions of the quasilinear equation involving the infinity Laplacian},  Nonlinear Analysis 74 (2011) 3478–3486.

\bibitem[K1]{K1} N. Katzourakis, \emph{$ L^\infty$-Variational Problems for Maps and the Aronsson PDE system}, J. Differential
Equations,Volume 253, Issue 7 (2012), 2123 - 2139.

\bibitem[K2]{K2} N. Katzourakis, \emph{An Introduction to Viscosity Solutions for Fully Nonlinear PDE with Applications to Calculus of Variations in $L^\infty$}, Springer Briefs in Mathematics,2015, DOI 10.1007/978-3-319-12829-0.

\bibitem[K3]{K3} N. Katzourakis, \emph{$\infty$-Minimal Submanifolds}, Proceedings of the AMS 142, 2797 - 2811 (2014).

\bibitem[K4]{K4} N. Katzourakis, \emph{On the Structure of $\infty$-Harmonic Maps}, Communications in PDE, Volume 39, Issue 11 (2014), 2091 - 2124.

\bibitem[K5]{K5} N. Katzourakis, \emph{Explicit 2D $\infty$-Harmonic Maps whose Interfaces have Junctions and Corners}, Comptes Rendus Acad. Sci. Paris, Ser.I, 351 (2013) 677 - 680.

\bibitem[K6]{K6} N. Katzourakis, \emph{Optimal $\infty$-Quasiconformal Immersions}, ESAIM Control Optim. Calc. Var.21 (2015), no. 2, 561 - 582.

\bibitem[K7]{K7} N. Katzourakis, \emph{ Nonuniqueness in Vector-valued Calculus of Variations in $L^\infty$ and some Linear Elliptic Systems}, Comm. on Pure and Appl. Anal., Vol. 14, 1, 313 - 327 (2015).

\bibitem[K8]{K8} N. Katzourakis, \emph{Generalised solutions for fully nonlinear PDE systems and existence uniqueness theorems}, Journal of Differential Equations 23 (2017), 641-686, DOI: 10.1016/j.jde.2017.02.048.

\bibitem[K9]{K9} N. Katzourakis, \emph{Absolutely minimising generalised solutions to the equations of vectorial Calculus of Variations in $L^\infty$}, Calculus of Variations and PDE 56 (1), 1 - 25 (2017) (DOI:10.1007/s00526-016-1099-z).

\bibitem[K10]{K10} N. Katzourakis, \emph{ A new characterisation of $\infty$-harmonic and p-harmonic maps via affine variations in $L^\infty$}, Electronic J. Differential Equations, Vol. 2017 (2017), No. 29, 1 - 19.

\bibitem[K11]{K11} N. Katzourakis, \emph{Solutions of vectorial Hamilton-Jacobi equations are rank-one Absolute Minimisers in $L^\infty$}, Advances in Nonlinear Analysis, in press.

\bibitem[KM]{KM} N. Katzourakis, J. Manfredi, \emph{Remarks on the Validity of the Maximum Principle for the $\infty$-Laplacian},  Le Matematiche, Vol. LXXI (2016) Fasc. I, 63 - 74, DOI: 10.4418/2016.71.1.5.

\bibitem[KMo]{KMo} N. Katzourakis, R. Moser, \emph{Existence, uniqueness and structure of second order absolute minimisers}, ArXiv preprint, \url{http://arxiv.org/pdf/1701.03348.pdf}.

\bibitem[KPa]{KPa} N. Katzourakis, E. Parini, \emph{The eigenvalue problem for the $\infty$-Bilaplacian}, ArXiv preprint, \url{https://arxiv.org/abs/1703.03648}.

\bibitem[KP1]{KP1} N. Katzourakis, T. Pryer, \emph{On the numerical approximation of $\infty$-Harmonic mappings}, Nonlinear Differential Equations \& Applications 23 (6), 1-23 (2016).

\bibitem[KP2]{KP2} N. Katzourakis, T. Pryer, \emph{2nd order $L^\infty$ variational problems and the $\infty$-Polylaplacian}, ArXiv preprint, \url{http://arxiv.org/pdf/1605.07880.pdf}.

\bibitem[KP3]{KP3} N. Katzourakis, T. Pryer, \emph{On the numerical approximation of $\infty$-Biharmonic and p-Biharmonic functions}, ArXiv preprint, \url{https://arxiv.org/pdf/1701.07415.pdf}.

\bibitem[VHV]{VHV} M.-F. Bidaut-Véron, M. Garcia-Huidobro, L. Véron, \emph{Separable $\infty$-Harmonic functions in cones},  ArXiv preprint, \url{https://arxiv.org/pdf/1703.07297.pdf}.
 
\bibitem[S]{S} O.Savin, \emph{$C^{1}$ regularity for $\infty$-Harmonic functions in two dimensions}, Archive for Rational Mech and Analysis 176 (2005) 351-361.

\end{thebibliography}

\end{document}